\documentclass[12pt,leqno]{amsart}

\usepackage{float}
\usepackage{amssymb,amsmath,amsthm,amscd}
\usepackage{mathrsfs}

\usepackage[all]{xy} 
\CompileMatrices

\numberwithin{equation}{section}

\newcommand{\Q}{\mathbb {Q}}
\newcommand{\BQ}{{\mathbb Q}}
\newcommand{\Z}{{\mathbb Z}}

\newcommand{\RO}{\mathscr{O}_{{\mathrm{int}}}}

\newcommand{\Gt}{{\mathfrak{t}}}
\newcommand{\td}{\Gt^*}

\newcommand{\seteq}{\mathbin{:=}}

\theoremstyle{plain}
\newtheorem{lemma}{Lemma}[section]
\newtheorem{prop}[lemma]{Proposition}
\newtheorem{theorem}[lemma]{Theorem}
\newcommand{\Prop}{\begin{prop}}
\newcommand{\enprop}{\end{prop}}
\newcommand{\Lemma}{\begin{lemma}}
\newcommand{\enlemma}{\end{lemma}}
\newcommand{\Th}{\begin{theorem}}
\newcommand{\enth}{\end{theorem}}
\newtheorem{corollary}[lemma]{Corollary}
\newcommand{\Cor}{\begin{corollary}}
\newcommand{\encor}{\end{corollary}}
\newtheorem{definition}[lemma]{Definition}

\newcommand{\Def}{\begin{definition}}
\newcommand{\edf}{\end{definition}}

\theoremstyle{definition}
\newtheorem{remark}[lemma]{Remark}

\newtheorem{conjecture}[lemma]{Conjecture}

\newtheorem*{thm}{Theorem \ref{th:main}}

\newcommand{\g}{{\mathfrak{g}}}
\newcommand{\Gg}{{\mathfrak{g}}}

\newcommand{\isoto}[1][]%
{{\mathop{\buildrel{\sim}\over\longrightarrow}\limits_{#1}}}
\renewcommand{\hom}{\operatorname{\it \mathscr{H}\kern-.25em om}}

\newcommand{\eq}{\begin{eqnarray}}
\newcommand{\eneq}{\end{eqnarray}}

\newcommand{\eqn}{\begin{eqnarray*}}
\newcommand{\eneqn}{\end{eqnarray*}}

\newenvironment{tenumerate}{
  \begin{enumerate}
  
  }{\end{enumerate}}

\newcommand{\bnum}{\begin{tenumerate}}
\newcommand{\enum}{\end{tenumerate}}
\newenvironment{anumerate}{
  \begin{enumerate}
  
  }{\end{enumerate}}

\newcommand{\bia}{\begin{anumerate}}
\newcommand{\eia}{\end{anumerate}}
\newcommand{\on}{\operatorname}
\newcommand{\Ker}{\on{Ker}}

\newcommand{\bni}{\begin{tenumerate}}
\newcommand{\eni}{\end{tenumerate}}

\newcommand{\QED}{\end{proof}}
\newcommand{\Proof}{\begin{proof}}

\renewcommand{\ss}{\smallskip}

\newcommand{\cl}{\colon}
\newcommand{\To}[1][\phantom{aaaa}]{\xrightarrow{\,#1\,}}

\newcommand{\ba}{\begin{array}}
\newcommand{\ea}{\end{array}}

\renewcommand{\phi}{\varphi}
\newcommand{\bi}{\bni}
\newcommand{\ei}{\eni}

\newcommand{\epi}{\twoheadrightarrow}

\newcommand{\set}[2]{\left\{#1 \mathbin; #2 \right\}}

\newcommand{\eqsub}{\begin{subequations}\begin{eqnarray}}
\newcommand{\eneqsub}{\end{eqnarray}\end{subequations}}

\newcommand{\ol}{\overline}

\newcommand{\A}{\mathbf{A}}

\renewcommand{\le}{\leqslant}
\renewcommand{\ge}{\geqslant}

\newcommand{\nc}{\newcommand}
\nc{\la}{\lambda}
\nc{\lam}{\lambda}
\nc{\U}{U_q(\g)}
\nc{\te}{\tilde{e}}
\nc{\tei}{\tilde{e}_i}
\nc{\tf}{\tilde{f}}
\nc{\tfi}{\tilde{f}_i}
\nc{\tU}{\widetilde U_q(\g)}

\nc{\BZ}{{\mathbb{Z}}}
\nc{\al}{\alpha}
\nc{\qs}{{q_{\mathrm{s}}}}
\nc{\lan}{\langle}
\nc{\ran}{\rangle}
\nc{\re}{{\mathrm{re}}}
\nc{\wt}{\operatorname{wt}}
\nc{\Uf}{U^-_q(\g)}
\nc{\Ue}{U^+_q(\g)}
\nc{\eps}{\varepsilon}
\nc{\vphi}{\varphi}
\nc{\sphi}{\varphi^*}
\nc{\seps}{\varepsilon^*}

\nc{\nn}{\nonumber}
\def\max{{\mathop{\mathrm{max}}}}
\nc{\vp}{\varpi}
\renewcommand{\L}{{\Lambda_0}}
\nc{\cls}{{\operatorname{cl}}}
\renewcommand{\O}{\operatorname{O}}
\nc{\Wt}{{\operatorname{Wt}}}
\nc{\Us}{U'_q(\g)}
\nc{\La}{\Lambda}
\nc{\ro}{{\rm(}}
\nc{\rf}{{\rm)}}
\nc{\norm}{{\mathrm{norm}}}
\nc{\qbox}{\quad\mbox}
\nc{\braid}{{\mathfrak{B}}}
\nc{\Ad}{\operatorname{Ad}}
\nc{\Aut}{\operatorname{Aut}}
\nc{\dt}[1]{\tilde{\tilde #1}}
\nc{\Sn}{S^{{\mathrm{norm}}}}
\nc{\aff}{{\mathrm{aff}}}
\nc{\rk}{{\mathrm{rk}}}
\nc{\tQ}{\widetilde{Q}}
\nc{\tP}{\widetilde{P}}
\nc{\tW}{\widetilde{W}}
\nc{\Dyn}{\mathrm{Dyn}}
\nc{\tD}{\widetilde{\Delta}}
\nc{\height}{{\operatorname{ht}}}
\nc{\bl}{\bigl}
\nc{\br}{\bigr}

\begin{document}

\title[Fundamental representations and Demazure modules]%
{Level zero fundamental representations over quantized affine algebras
and Demazure modules}
\author{Masaki KASHIWARA}
\address{Research Institute for Mathematical Sciences,
Kyoto University, Kyoto 606, Japan
}
\thanks{This research is partially supported by 
Grant-in-Aid for Scientific Research (B1)13440006,
Japan Society for the Promotion of Science.}

\begin{abstract}
Let $W(\vp_k)$ be the finite-dimensional
irreducible module over a quantized affine algebra $\Us$
with the fundamental weight $\vp_k$ as an extremal weight.
We show that its crystal
$B(W(\vp_k))$ is isomorphic to the
Demazure crystal $B^-(-\L+\vp_k)$.
This is derived from the following general result:
for a dominant integral weight $\la$ 
and an integral weight $\mu$,
there exists a unique homomorphism
$\U (u_\la\otimes u_\mu)\to V(\la+\mu)$
that sends $u_\la\otimes u_\mu$ to $u_{\la+\mu}$.
Here $V(\la)$ is the extremal weight module with
$\la$ as an extremal weight, 
and $u_\la\in V(\la)$ is the extremal weight vector of
weight $\la$.
\end{abstract}

\maketitle

\section{introduction}
The finite-dimensional representations of
quantized affine algebras $\Us$ 
are extensively studied in connection with exactly solvable models.
It is expected that there exists
a ``good'' finite-dimensional 
$\Us$-module $W(m\vp_k)$ with a multiple $m\vp_k$  of 
a fundamental weight $\vp_k$
as an extremal weight. This module
is good in the sense that it is irreducible and 
it has a crystal base
and moreover a global basis.

In the untwisted case,
its conjectural character formula is given by
Kirillov--Reshetikhin (\cite{KR}, see also \cite{kelber}), and
its conjectural fusion construction is given by Kuniba--Nakanishi--Suzuki
(\cite{KNS}).
It is proved by Nakajima (\cite{N2}) that the fusion construction 
gives irreducible modules
with the expected character in the simply laced case,
and by Chari (\cite{Chari}) in some cases.

It is also expected that any ``good'' finite-dimensional $\Us$-module
is a tensor product of modules of the above type.

It is also conjectured in \cite{HKOTY, HKOTT} 
that the $\Us$-modules $W(m\vp_k)$
has a perfect crystal of level $\ell$
if and only if $m=\ell c^\vee_k$ ( $c^\vee_k\seteq\max(1,2/(\al_k,\al_k))$).
Moreover it is conjectured that the crystal base
$B(W(\ell c^\vee_k\vp_k))$ is isomorphic to the
Demazure crystal $B^-(-\ell\L+\ell c^\vee_k\vp_k)$ if we forget the $0$-arrows.
Here, for an integral weight $\la$, $B^\pm(\la)$ denotes the crystal for
the $U_q^{\pm}(\Gg)$-module
generated by the extremal vector
with weight $\la$.
They are proved in certain cases (\cite{(KMN)^2,(KMN)^2_2}).
More general relations of perfect crystals and Demazure crystals
are discussed in \cite{KMOTU}.

In this paper we show that
$B(W(\vp_k))$ is isomorphic to the
Demazure crystal $B^-(-\L+\vp_k)$,
or equivalently 
$B(W(-\vp_k))$ is isomorphic to the
Demazure crystal $B^+(\L-\vp_k)$ (Corollary \ref{cor:fund}).

The main ingredient is the following theorem,
which the author started to study in order to answer
a question raised by Miwa et al:

\begin{thm}
Let $\U$ be a quantized affine algebra.
Let $\la\in P^+$ be a dominant integral weight 
and $\mu\in P$ an integral weight.
Then there exists a unique homomorphism
$V(\la)\otimes V(\mu)\supset\U (u_\la\otimes u_\mu)\To V(\la+\mu)$
that sends $u_\la\otimes u_\mu$ to $u_{\la+\mu}$.
Moreover this morphism is compatible with global bases.
\end{thm}
Here $V(\la)$ is the extremal weight module with
$\la$ as an extremal weight, 
and $u_\la\in V(\la)$ is the extremal weight vector of
weight $\la$.

\bigskip
\noindent
{\em Acknowledgment\quad}
The author thanks B. Feigin, M. Jimbo, T. Miwa, E. Mukhin, Y. Takeyama, 
and M. Okado
for helpful discussions. 

\section{Review on crystal bases and global bases}\label{sec:rev1}

In this section, we shall review briefly 
the quantized universal enveloping algebras and crystal bases.
We refer the reader to \cite{banff,quant,modified,KF,cours,GL}.

\subsection{Quantized universal enveloping algebras}
We shall define the quantized universal enveloping algebra
$U_q(\g)$.
Assume that we are given the following data.
\begin{eqnarray*}
&&P:\text{a free $\BZ$-module (called a weight lattice),}\\
&&I:\text{an index set (for simple roots),}\\
&&\al_i\in P\ \text{for}\ i\in I\ \text{(called a simple root),}\\
&&h_i\in P^*\seteq{\mathrm{Hom}}_\BZ(P,\BZ)\ \text{(called a simple coroot),}\\
&&(\,\cdot\,,\,\cdot\,)\colon P\times P\to \BQ\quad \text{a bilinear symmetric form.}
\end{eqnarray*}
We shall denote by $\langle \,\cdot\,,\,\cdot\,\rangle 
\colon  P^*\times P\to \BZ$ the
canonical pairing. 

The data above are assumed to satisfy the following axioms.
\begin{eqnarray}
&&\ba{l}
(\al_i,\al_i)>0\quad\text{for any $i\in I$,}\\[2pt]
(\al_i,\al_j)\le 0\quad 
\text{for any $i$, $j\in I$ with $i\not= j$,}\\[2pt]
\langle h_i,\la\rangle =\dfrac{2(\al_i,\lam)}{(\al_i,\al_i)}
\quad\text{for any $i\in I$ and $\lam\in P$.}
\ea
\end{eqnarray}

Let us take a positive integer  $d$
such that $(\al_i,\al_i)/2\in\Z\,d^{-1}$ for any $i\in I$.
Now let $q$ be an indeterminate and set 
\eq
&&\text{$K=\Q(\qs)$ where
$\qs=q^{1/d}$.}
\eneq
We define its subrings $\A_0$, $\A_\infty$ and $\A$ as follows.
\eq
&&\ba{rcl}
\A_0&=&\set{f/g}{f,g\in \Q[\qs],\, g(0)\not=0},\\[3pt]
\A_\infty&=&\set{f/g}{f,g\in \Q[\qs^{-1}],\, g(0)\not=0},\\[3pt]
\A&=&\Q[\qs,\qs^{-1}].
\ea
\eneq

\begin{definition}\label{U_q(g)}
The quantized universal enveloping algebra $U_q(\g)$ is the algebra over
$K$ generated by the symbols $e_i,f_i\ (i\in I)$ and $q(h)\ (h\in d^{-1}P^*)$ 
with the following defining relations.
\begin{enumerate}

\item $q(h_1)q(h_2)=q(h_1+h_2)$ for $h_1,h_2\in d^{-1}P^*$,
and $q(h)=1$ for $h=0$.

\item
$q(h)e_i\,q(h)^{-1}=q^{\langle h,\al_i\rangle}\,e_i\ $
and $\ q(h)f_i\,q(h)^{-1}=q^{-\langle h,\al_i\rangle}f_i\ $
for any $i\in I$ and $h\in d^{-1}P^*$.

\item\label{even} $\lbrack e_i,f_j\rbrack
=\delta_{ij}\dfrac{t_i-t_i^{-1}}{q_i-q_i^{-1}}$
for $i$, $j\in I$. Here $q_i=q^{(\al_i,\al_i)/2}$ and
$t_i=q(\frac{(\al_i,\al_i)}{2}h_i)$.

\item {\rm(}{\em Serre relation}{\rm)} For $i\not= j$,
\begin{eqnarray*}
&&\sum^b_{k=0}(-1)^ke^{(k)}_ie_je^{(b-k)}_i=\sum^b_{k=0}(-1)^kf^{(k)}_i
f_jf_i^{(b-k)}=0.
\end{eqnarray*}
Here $b=1-\langle h_i,\al_j\rangle$ and
\begin{eqnarray*}
e^{(k)}_i=e^k_i/\lbrack k\rbrack_i!\ ,&& f^{(k)}_i=f^k_i/\lbrack k\rbrack_i!\ ,\\
\lbrack k\rbrack_i=(q^k_i-q^{-k}_i)/(q_i-q^{-1}_i)\ ,
&&\lbrack k\rbrack_i!=\lbrack 1\rbrack_i\cdots \lbrack k\rbrack_i\,.
\end{eqnarray*}
\end{enumerate}
\end{definition}

For $i\in I$, we denote by $\U_i$ the subalgebra of $\U$ generated by $e_i$,
$f_i$ and $q(h)$ ($h\in d^{-1}P^*$).

Let us denote by $W$ the Weyl group, the subgroup of $GL(P)$ generated
by the simple reflections $s_i$:
$s_i(\lam)=\lam-\lan h_i,\lam\ran\alpha_i$.

Let $\Delta\subset Q\seteq\sum_i\Z\alpha_i$ be the set of roots.
Let $\Delta^\pm\seteq\Delta\cap Q_\pm$ be the set of positive 
and negative roots, respectively.
Here
$Q_\pm\seteq\pm\sum_i\Z_{\ge0}\alpha_i$.
Let $\Delta^\re$ be the set of real roots, and set
$\Delta_\pm^\re\seteq\Delta_\pm\cap\Delta^\re$.

\subsection{Braid group action on integrable modules}
The $q$-analogue of the action of the Weyl group
is introduced in \cite{GL,S}.
We define a $q$-analog of the exponential function by
\eq
\exp_q(x)=\sum_{n=0}^\infty {\dfrac{q^{n(n-1)/2}x^n}{[n]!}}\,.
\eneq
This satisfies the following equations:
\eq&&
\ba{lll}
&&\exp_q(x)\exp_{q}(y)=\exp_q(x+y)\qbox{if $xy=q^2yx$,}
\label{eq:add}\\[3pt]
&&\exp_q(x)\exp_{q^{-1}}(y)=\sum_{n=0}^\infty
\dfrac{1}{[n]!}\prod_{\nu=0}^{n-1}(q^\nu x+q^{-\nu}y)
\qbox{if $[x,y]=0$,}\\[3pt]
&&\exp_q(x)\exp_{q^{-1}}(-x)=1,\\[5pt]
&&\exp_q(x)=\Big(1+(1-q^2)x\Big)\exp_q(q^{2}x),\\[5pt]
&&\exp_q(x)=\prod_{n=0}^\infty\Big(1+q^{2n}(1-q^2)x\Big)
\qbox{for $|q|<1$,}
\ea
\eneq

For $i\in I$,
we set
\eq
&&\ba{rcl}
\hspace*{3em}S_i&=&\exp_{q_i^{-1}}(q_i^{-1}e_it_i^{-1})
\exp_{q_i^{-1}}(-f_i)\exp_{q_i^{-1}}(q_ie_it_i)
\,q_i^{h_i(h_i+1)/2}\\
&=&\exp_{q_i^{-1}}(-q_i^{-1}f_it_i)
\exp_{q_i^{-1}}(e_i)\exp_{q_i^{-1}}(-q_if_it_i^{-1})
\,q_i^{h_i(h_i+1)/2}. 
\ea
\eneq
We regard $S_i$ 
as an endomorphism of 
integrable $U_q(\Gg)$-modules,
and $q_i^{h_i(h_i+1)/2}$ acts on the weight space
of weight $\lam$ by the multiplication of
$q_i^{\lan h_i,\lam\ran(\lan h_i,\lam\ran+1)/2}$.

On the $(l+1)$-dimensional irreducible representation of
$U_q(\Gg)_i$
with a highest weight vector $u_0^{(l)}$ 
and $u_k^{(l)}=f_i^{(k)}u_0^{(l)}$,
\eq
&&\ba{l}
S_i(u_k^{(l)})=(-1)^{l-k}q_i^{(l-k)(k+1)}u_{l-k}^{(l)},\label{eq:br:vec}%
\\[5pt]
\ea
\eneq
Hence, $S_i$ sends the weight space of weight $\lam$ 
to the weight space of weight $s_i\lam$.
By the above formula, we have
\eq\label{eq:2:1}
&&\mbox{$S_iu_l^{(l)}=u_0^{(l)}$\quad
and\quad $S_iu_0^{(l)}=(-q_i)^lu_l^{(l)}$.}
\eneq

Since $\{S_i\}$ satisfies the braid relations,
we can extend the actions of $S_i$ on
integrable modules
to the action of the braid group
by
\eqn
&&
\begin{matrix}
S_{ww'}&=&S_w\circ S_{w'} \qbox{if $l(ww')=l(w)+l(w')$,}\\
S_{s_i}&=&S_i\,.\hfill
\end{matrix}
\eneqn

\subsection{Braid group action on $\U$}

We define the ring automorphism $T_i$ of $U_q(\Gg)$
by
\eq
T_i(q)&=&q\\
T_i(q(h))&=&q(s_ih),\\
T_i(e_i)&=&-f_it_i,\\
T_i(f_i)&=&-t_i^{-1}e_i,\\
T_i(e_j)&=&\sum_{k=0}^{-\lan h_i,\alpha_j\ran}(-1)^kq_i^{-k}e_i^{(-\lan h_i,\alpha_j\ran-k)}
e_je_i^{(k)},\\
T_i(f_j)&=&\sum_{k=0}^{-\lan h_i,\alpha_j\ran}(-1)^kq_i^{k}f_i^{(k)}
f_jf_i^{(-\lan h_i,\alpha_j\ran-k)}\hbox{ for $i\neq j$.}
\eneq
Then it is well-defined, and it satisfies
\eq
&&T_i(P)u=S_iPS_i^{-1}u
\eneq
for any $P\in U_q(\Gg)$ and any element $u$
of an integrable $U_q(\Gg)$-module.

The operator $T_i$ is invertible and its inverse is given as follows.
\eq
T_i^{-1}(q(h))&=&q(s_ih),\\
T_i^{-1}(e_i)&=&-t_i^{-1}f_i,\\
T_i^{-1}(f_i)&=&-e_it_i,\\
T_i^{-1}(e_j)&=&\sum_{k=0}^{-\lan h_i,\alpha_j\ran}(-1)^kq_i^{-k}e_i^{(k)}
e_je_i^{(-\lan h_i,\alpha_j\ran-k)},\\
T_i^{-1}(f_j)&=&\sum_{k=0}^{-\lan h_i,\alpha_j\ran}(-1)^kq_i^{k}f_i^{(-\lan h_i,\alpha_j\ran-k)}
f_jf_i^{(k)}.
\eneq

We can extend the action $T_i$ to the action of the braid group
by
\eqn
&&
\begin{matrix}
T_{ww'}&=&T_w\circ T_{w'} \qbox{if $l(ww')=l(w)+l(w')$,}\\
T_{s_i}&=&T_i\,.\hfill
\end{matrix}
\eneqn
%

The following proposition is proved in \cite{GL}.

\Prop
For $w\in W$ and $i,j\in I$ such that
$w\alpha_i=\alpha_j$, we have
$$T_we_i=T_{w^{-1}}^{-1}e_i=e_j\qbox{and}\quad
T_wf_i=T_{w^{-1}}^{-1}f_i=f_j.$$
\enprop

\subsection{Crystals}
We shall not review the notion of crystals,
but refer the reader to
\cite{banff,quant,modified,cours}.
For a subset $J$ of $I$, let us denote by $U_q(\Gg_J)$ 
the subalgebra of $\U$ generated by
$e_j$, $f_j$ ($j\in J$) and $q(h)$ ($h\in d^{-1}P^*$).
We say that a crystal $B$ over $\U$
is a {\em regular crystal\/} 
if, for any $J{\subset}I$ of finite-dimensional type,
$B$ is, as a crystal over $U_q(\Gg_J)$, 
isomorphic to a  crystal base associated with
an integrable $U_q(\Gg_J)$-module.

By \cite{modified},
the Weyl group $W$ acts on any regular crystal.
This action $S$ is given by
\eqn
&&S_{s_i}b=
\begin{cases}
\tf_i^{\lan h_i,\wt(b)\ran}b
&\mbox{if $\lan h_i,\wt(b)\ran\ge 0$,}\\
\te_i^{-\lan h_i,\wt(b)\ran}b
&\mbox{if $\lan h_i,\wt(b)\ran\le 0$.}
\end{cases}
\eneqn

Let us denote by
$\Uf$ (resp.\ $\Ue$) 
the subalgebra of $\U$ generated by the $f_i$'s
(resp.\ by the $e_i$'s).
Then $\Uf$ has a crystal base denoted by
$B(\infty)$ (\cite{quant}). A unique vector of
$B(\infty)$ with weight $0$ is denoted by $u_\infty$.
Similarly $\Ue$ has a crystal base denoted by
$B(-\infty)$, and a unique vector of
$B(-\infty)$ with weight $0$ is denoted by $u_{-\infty}$.

Let $\psi$ be the ring automorphism of $\U$ that sends
$\qs$, $e_i$, $f_i$ and $q(h)$ to $\qs$, $f_i$, $e_i$ and $q(-h)$.
It induces bijections $\Uf\isoto \Ue$ and
$B(\infty)\isoto B(-\infty)$ by which
$u_\infty$, $\te_i$, $\tf_i$, $\eps_i$, $\vphi_i$, $\wt$
correspond to $u_{-\infty}$, $\tf_i$, $\te_i$, $\vphi_i$, $\eps_i$, 
$-\wt$.

Let $\tU$ be the modified quantized
universal enveloping algebra $\oplus_{\lam\in P}\U a_\la$
(see \cite{modified}).
The elements $a_\la$, the projectors to the weight $\la$-space,
satisfy $a_\la\cdot a_\mu=\delta_{\la,\mu}a_\la$ and
$a_\la P=P a_{\la-\wt(P)}$ for $P\in \U$.

Then $\tU$ has a crystal base $(L(\tU),B(\tU))$.
As a crystal,
$B(\tU)$ is regular and isomorphic to
$$\bigsqcup_{\lam\in P}B(\infty)\otimes T_\lam\otimes B(-\infty).$$
Here, $T_\lam$ is the crystal consisting of a single element
$t_\lam$ with $\eps_i(t_\lam)=\vphi_i(t_\lam)=-\infty$ 
and $\wt(t_\lam)=\lam$.

Let $*$ be the anti-involution of $\U$ that sends
$q(h)$ to  $q(-h)$, and $\qs$, $e_i$, $f_i$ to themselves.
The involution $*$ of $\U$ induces
an involution $*$ on $B(\infty)$, $B(-\infty)$, $B(\tU)$.
Then $\te_i^*=*\circ\te_i\circ*$, etc.
give another crystal structure on $B(\infty)$, $B(-\infty)$, $B(\tU)$.
We call it {\em the star crystal structure}.
These two crystal structures on $B(\tU)$ are compatible,
and  $B(\tU)$ may be considered as a
crystal over $\g\oplus \g$,
which corresponds to the $\U$-bimodule structure on $\tU$.
Hence, for example,
$S^*_w$, the Weyl group action on $B(\tU)$
with respect to the star crystal structure
is a crystal automorphism of $B(\tU)$ with respect
to the original crystal structure.
In particular, the two Weyl group actions 
$S_w$ and $S_{w'}^*$ commute with each other.

\subsection{Global bases}\label{subsec:global}
Recall that $\A_0\subset K$ is the subring of $K$ consisting of rational
functions in $\qs$ without pole at $\qs=0$.
Let $-$ be the automorphism of $K$ sending $\qs$ to $\qs^{-1}$.
Then $\ol{\A_0}$ coincides with the ring $\A_\infty$ of rational functions 
regular at $\qs=\infty$.
Set $\A\seteq\Q[\qs,\qs^{-1}]$.
Let $V$ be a vector space over $K$,
$L_0$ an $A$-submodule of $V$,
$L_\infty$ an $\A_\infty$- submodule, and
$V_\A$ a $\A$-submodule.
Set $E:=L_0\cap L_\infty\cap V_\A$.

\Def[\cite{quant}]
We say that $(L_0,L_\infty,V_\A)$ is {\em balanced}
if each of $L_0$, $L_\infty$ and $V_\A$
generates $V$ as a $K$-vector space,
and if one of the following equivalent conditions is satisfied.
\bnum
\item
$E \to L_0/\qs L_0$ is an isomorphism.
\item
$E \to L_\infty/\qs^{-1}L_\infty$ is an isomorphism.
\item
$(L_0\cap V_\A)\oplus
(\qs^{-1} L_\infty \cap V_\A) \to V_\A$
      is an isomorphism.
\item
$\A_0\otimes_\Q E \to L_0$, $\A_\infty\otimes_\Q E \to L_\infty$,
        $\A\otimes_\Q E \to V_\A$ and $K \otimes_\Q E \to V$
are isomorphisms.
\enum
\edf

Let $-$ be the ring automorphism of $\U$ sending
$\qs$, $q(h)$, $e_i$, $f_i$ to $\qs^{-1}$, $q(-h)$, $e_i$, $f_i$.

Let $\U_\A$ be the $\A$-subalgebra of
$\U$ generated by $e_i^{(n)}$, $f_i^{(n)}$
and $q(h)$ ($h\in d^{-1}P^*$). 

Let $M$ be a $\U$-module.
Let $-$ be an involution of $M$ satisfying $(au)^-=\bar a\bar u$
for any $a\in\U$ and $u\in M$.
We call in this paper such an involution a {\em bar involution}.
Let $(L(M),B(M))$ be a crystal base of an integrable $\U$-module $M$.

Let $M_\A$ be a $\U_\A$-submodule of $M$
such that 
\eq\label{eq:zform}
&&
\text{$(M_\A){}^-=M_\A$, 
and $(u-\ol{u})\in (\qs-1)M_\A$ for every $u\in M_\A$.}
\eneq

\Def
A $\U$-module $M$ endowed 
with $(L(M),B(M),M_\A,-)$ as above is called
with a {\em global basis},
if $(L(M),L(M)^-,M_\A)$ is balanced, 
\edf
In such a case, let $G\colon L(M)/\qs L(M)\isoto E:=L(M)\cap L(M)^-
\cap M_\A$ 
be the inverse of $E\isoto L(M)/\qs L(M)$.
Then $\{G(b);b\in B(M)\}$ forms a basis of $M$.
We call this basis a (lower) {\em global basis}.
The global basis enjoys the following properties
(see \cite{quant,global}):
\bnum
\item $\ol{G(b)}=G(b)$ for any $b\in B(M)$.
\item
For any $n\in\Z_{\ge0}$,
$\{G(b);\eps_i(b)\ge n\}$ is a basis
of the $\A$-submodule $\sum_{m\ge n}f_i^{(m)}M_\A$.
\item
for any $i\in I$ and $b\in B(M)$,
we have
\[f_iG(b)=[1+\eps_i(b)]_iG(\tf_ib)+\sum_{b'}F^{i}_{b,b'}G(b').\]
Here the sum ranges over
$b'\in B(M)$ such that $\eps_i(b')>1+\eps_i(b)$.
The coefficient $F^i_{b,b'}$ belongs to
$\qs q_i^{1-\eps_i(b')}\Q[\qs]$.
Similarly for $e_iG(b)$.
\enum

\smallskip
Let $M$ and $N$ be $\U$-modules with global bases.
We say that a $\U$-morphism $f\cl M\to N$ is 
{\em compatible with global bases}
if it satisfies the following conditions:
\bi
\item
If $u$ is a global basis vector of $M$, then $f(u)$ is a global basis vector
of $N$ or $0$.
\item
If a pair of global basis vectors $u$ and $v$ of $M$ satisfies
$f(u)=f(v)\not=0$, then $u=v$.
\ei
These conditions are equivalent to the following set of conditions:
\bia
\item
$f$ commutes with the bar involutions.
\item
$f$ sends $L(M)$ to $L(N)$ and $M_\A$ to $N_\A$.
\item
The induced morphism $\ol{f}\cl L(M)/\qs L(M)\to L(N)/\qs L(N)$
sends $B(M)$ to $B(N)\cup\{0\}$.
\item
$\Ker(f)$ is generated by a part of the global basis of $M$.
\eia
In such a case, $f(M)$ 
has a global basis,
and we have
$$B(M)\supset B(f(M))\subset B(N).$$
If $f$ is a monomorphism then $B(M)\simeq B(f(M))\subset B(N)$,
and if $f$ is an epimorphism then $B(M)\supset B(f(M))\simeq B(N)$.

\subsection{Extremal vectors}
Let $M$ be an integrable $\U$-module.
A non-zero vector $u\in M$ of weight $\lambda\in P$
is called {\em extremal} (see \cite{modified}),
if we can find
a subset $F$ of non-zero weight vectors in $M$ containing $u$ 
and satisfying the following properties:
\eq
&&\ba{l}
\hbox{if $v\in F$ and $i$ satisfy $\lan h_i,\wt(v)\ran\ge 0$, then
$e_i v=0$ and $f_i^{(\lan h_i,\wt(v)\ran)}v\in F$,}\\
\hbox{if $v\in F$ and $i$ satisfy $\lan h_i,\wt(v)\ran\le 0$, then
$f_i v=0$ and $e_i^{(-\lan h_i,\wt(v)\ran)}v\in F$,}
\ea
\eneq
The Weyl group $W$ acts on the set of extremal vectors by
\eq
&&\ba{l}
\hbox{if $\lan h_i,\wt(u)\ran\ge 0$, then
$\Sn_{s_i} u=f_i^{(\lan h_i,\wt(u))}u$,}\\
\hbox{if $\lan h_i,\wt(u)
\ran\le 0$, then
$\Sn_{s_i}u=e_i^{(-\lan h_i,\wt(u)\ran)}u$.}
\ea
\eneq
We have $\wt(\Sn_wu)=w\wt(u)$ for $w\in W$.
Note that, by \eqref{eq:br:vec},
$\Sn_w u$ is equal to $S_w u$ up to a non-zero constant multiple.

Similarly, for a vector $b$ of a regular crystal $B$ with weight $\lam$,
we say that $b$
is an extremal vector
if it satisfies the following similar conditions:
\eq
&&\ba{l}
\hbox{if $w\in W$ and $i\in I$ satisfy $\lan h_i,w\lam\ran\ge 0$, then
$\te_iS_w b=0$,}\\
\hbox{if $w\in W$ and $i\in I$ satisfy $\lan h_i,w\lam\ran\le 0$ then
$\tf_iS_w b=0$.}
\ea
\eneq

For $\lam\in P$,
let us denote by $V(\lam)$
the $\U$-module
generated by $u_\lam$
with the defining relation that
$u_\lam$ is an extremal vector of weight $\lam$.
This is in fact infinitely many linear relations on $u_\lam$.

For a dominant weight $\lam$,
$V(\lam)$ is an irreducible highest weight module 
with highest weight $\lam$,
and $V(-\lam)$ is an irreducible lowest weight module
with lowest weight $-\lam$.

We proved in \cite{modified}
\footnote{ In \cite{modified}, it is denoted by $V^\max(\lam)$,
because I thought there would be a natural $\U$-module 
whose crystal base is the connected component of $B(\lam)$.}
that $V(\lam)$ has a global basis
$(L(\lam),B(\lam))$.
We denote by the same letter $u_\lam$ the element of $B(\lam)$
corresponding to $u_\lam\in V(\lam)$.
Moreover $\U a_\la\to V(\la)$ ($a_\la\mapsto u_\la$)
is compatible with global bases.
Hence the crystal $B(\lam)$ is isomorphic to
the subcrystal of
$B(\infty)\otimes t_\lam\otimes B(-\infty)$
consisting of vectors $b$ such that $b^*$ is an extremal vector
of weight $-\lam$.
By this embedding, $u_\lam\in B(\lam)$ corresponds to 
$u_\infty\otimes t_\lam\otimes u_{-\infty}$.

Note that 
\eq\label{eq:gbu}
&&
\Ue u_\la=
\bigoplus_{b\in B(\la)\cap \bl(u_\infty\otimes t_\lam\otimes B(-\infty)\br)}
KG(b).
\eneq

For any $w\in W$, $u_\lam\mapsto S^\norm_{w^{-1}}u_{w\lam}$
gives an isomorphism of $\U$-modules:
\eqn
&&V(\lam)\isoto V(w\lam).
\eneqn
This is compatible with global bases.
Similarly, 
letting $S^*_w$ be the Weyl group action on $B(\tU)$
with respect to the  star crystal structure
and regarding $B(\lam)$ as a subcrystal of $B(\tU)$,
$S^*_w\colon B(\tU)\isoto B(\tU)$ induces an isomorphism of crystals
\eq\label{eq:weyl}
S^*_w\colon B(\lam)\isoto B(w\lam).
\eneq
This coincides with the crystal isomorphism induced
by $V(\la)\isoto V(w\la)$.
Note that we have
$$S_wS_w^*(u_\infty\otimes t_\la\otimes u_{-\infty})
=u_\infty\otimes t_{w\la}\otimes u_{-\infty}.$$

\subsection{Global bases of tensor products}\label{sec:gltens}
Let us recall the following results proved by Lusztig (\cite{GL}).
Let $\RO$ be the category of integrable $\U$-modules which are a direct sum of
$V(\la)$'s ($\la\in P^+$).
Similarly let $\RO^-$ be the category of integrable $\U$-modules 
which are a direct sum of $V(\la)$'s ($\la\in P^-$).
Let $M$ and $N$ be $\U$-modules.
Assume that $M$ and $N$ have bar involutions,
and that either $M\in\RO$ or $N\in \RO^-$.
Then there exists a unique bar involution on $M\otimes N$ such that
$$(u\otimes v)^-=\bar u\otimes \bar v\qquad\parbox[t]{300pt}{for every 
$u\in M$ and $v\in N$ such that
either $u$ is a highest weight vector or $v$ is a lowest weight vector.}$$

Assume further that $M$ and $N$ have a global basis.
Then $M\otimes N$ has a crystal base
$(L(M\otimes N),B(M\otimes N))\seteq
(L(M)\otimes_{\A_0} L(N), B(M)\otimes B(N))$,
and an $\A$-form $(M\otimes N)_\A=M_\A\otimes_\A N_\A$.
Then $M\otimes N$ has a global basis;
namely $(L(M\otimes N),L(M\otimes N)^-,(M\otimes N)_\A)$
is balanced. In particular, $V(\la)\otimes V(\mu)$ has a global basis
either if $\la$ is dominant or if $-\mu$ is dominant.

\bigskip
Let $\la\in P$.
Then for any pair of
dominant integral weights $\xi$ and $\eta$ such that $\la=\xi-\eta$,
$\U a_{\la}\to V(\xi)\otimes V(-\eta)$
($a_{\la}\mapsto u_\xi\otimes u_{-\eta}$)
is compatible with global bases.
Conversely the global basis of $\U a_{\la}$ is characterized 
by the above property.

\Lemma
For $\la\in P^+$ and $\mu\in P$,
\eq\label{mor:uvu} 
\U a_{\la+\mu}\to V(\la)\otimes \U a_\mu\quad
\ro a_{\la+\mu}\mapsto u_\la\otimes a_{\mu}\rf
\eneq
is compatible with global bases.
\enlemma
\proof
For dominant integral weights $\xi$ and $\eta$ such that $\mu=\xi-\eta$,
we have a diagram of morphisms compatible with crystal basses
except the dotted arrow:
$$
\xymatrix{
\U a_{\la+\mu}\ar@{>>}[r]
\ar@{.>}[d]&{V(\la+\xi)\otimes  V(-\eta)}\ar@{^{(}->}[d]\\
{V(\la)\otimes \U a_\mu}\ar@{>>}[r]&{V(\la)\otimes V(\xi)\otimes  V(-\eta)}
}$$
Hence the dotted arrow is compatible with crystal bases.
\qed

\smallskip
This morphism \eqref{mor:uvu} induces an embedding of crystals
$$\mbox{$B(\U a_{\la+\mu})\hookrightarrow B(\la)\otimes B(\U a_\mu)$
for $\la\in P^{+}$ and $\mu\in P$.}
$$
There exists an embedding $B(\infty)\hookrightarrow 
B(\la)\otimes B(\infty)\otimes T_{-\la}$,
and the above morphism coincides with the composition
\eqn
B(\U a_{\la+\mu})&\simeq&B(\infty)\otimes T_{\la+\mu}\otimes B(-\infty)
\hookrightarrow B(\la)\otimes B(\infty)\otimes T_{-\la}
\otimes T_{\la+\mu}\otimes B(-\infty)\\
&\simeq& B(\la)\otimes B(\infty)\otimes T_\mu\otimes B(-\infty)
\simeq B(\la)\otimes B(\U a_\mu).
\eneqn

\subsection{Demazure modules}
Let $M$ be an integrable $\U$-module with a global basis
$(L(M)$,$B(M)$,$M_\A,-)$.
Let $N$ be a $\Ue$-submodule of $M$.
We say that
$N$ is {\em compatible with the global basis of} $M$
if  there exists a subset $B(N)$ of $B(M)$ such that
$N=\oplus_{b\in B(N)} K G(b)$.

It is shown in \cite{KL} that 
\eq&&
\parbox{350pt}{$\tei B(N)\subset B(N)\cup\{0\}$, and
$\U N=\Uf N$ is also compatible with the global basis.}
\eneq
Namely there exists a subset $B(\U N)$ of $B(M)$ such that
$$\U N=\bigoplus_{b\in B(\U N)} K G(b).$$
Moreover we have
$$B(\U N)=\set{\tf_{i_1}\cdots\tf_{i_m}b}%
{m\ge 0, i_1,\ldots i_m\in I, b\in B(N)}\setminus\{0\}.$$

For $\la\in P$, the $U_q^{\pm}(\Gg)$-submodule $U_q^{\pm}(\Gg)u_\la$ 
of $V(\la)$
is compatible with the global basis of $V(\la)$ (see \eqref{eq:gbu}).

We set
$$B^{\pm}(\la)=B(U^{\pm}(\g)u_\la).$$
Regarding $B(\la)$ as a subset of 
$B(\U a_\la)=B(\infty)\otimes t_\la\otimes B(-\infty)$,
we have
$$B^+(\la)
=B(\la)\cap\bigl(u_\infty\otimes t_\la\otimes B(-\infty)\bigr)
\quad\mbox{and}\quad
B^-(\la)=B(\la)\cap\bigl(B(\infty)\otimes t_\la\otimes u_{-\infty}\bigr).$$

The subset $B^+(\la)$ satisfies the following properties.
\Lemma\label{lem:str}
\bi
\item
$\tei B^+(\la)\subset B^+(\la)\cup\{0\}$.
\item
For any $b\in B^+(\la)$,
if $\eps_i(b)>0$, then $\tfi b\in B^+(\la)\cup\{0\}$.
Or equivalently, for any $i$-string $S$ of $B(\la)$,
$S\cap B^+(\la)$ is either $S$ itself, the empty set or
the set consisting of the highest weight vector of $S$.
Here an {\em $i$-string} is a connected component 
with respect to the crystal structure over $\U_i$.
\ei
\enlemma

This is a consequence of
the following lemma. Note that 
$B(\Ue a_\la)=u_\infty\otimes T_\la\otimes B(-\infty)$.
\Lemma
\bi
\item
$\tei B(\Ue a_\la)\subset B(\Ue a_\la)\cup\{0\}$.
\item
For any $b\in B(\Ue a_\la)$,
if $\eps_i(b)>0$, then $\tfi b\in B(\Ue a_\la)\cup\{0\}$.
Or equivalently, for any $i$-string $S$ of $B(\U a_\la)$,
$S\cap B(\Ue a_\la)$ is either $S$ itself, the empty set or
the set consisting of the highest weight vector of $S$.
\ei
\enlemma

\proof
The first property is evident.
In order to prove (ii), write
$b=u_\infty\otimes t_\la\otimes b'$ with $b'\in B(-\infty)$.
Then $\eps_i(b)=\max(0,\eps_i(t_\la\otimes b'))$, and hence
$0=\phi_i(u_\infty)<\eps_i(t_\la\otimes b')$.
We have therefore
$\tfi b=u_\infty\otimes t_\la\otimes \tfi b'$.
\qed

\ss
Similar results hold for $B^-(\la)$ and $B(\Uf a_\la)$.

\Prop\label{prop:beta}
For $\beta\in\Delta^\re_+$ and $\lam\in P$,
assume $(\beta,\lam)\ge0$.
Then we have
$$\mbox{$S_{s_\beta}u_\lam\in \Uf u_\lam$\quad and \quad
$S_{s_\beta}u_\la\in B^-(\la)$.}$$
\enprop

\proof
We shall argue by the induction of $\height(\beta)$.
Let us take $i\in I$ such that
$\lan h_i,\beta\ran>0$.
If $\beta=\alpha_i$ then the assertion is trivial.
Otherwise we have $\gamma\seteq s_i(\beta)\in\Delta^\re_+$.
Since $(\gamma,s_i\lam)=(\beta,\lam)\ge0$,
the induction hypothesis implies that
\eq\label{eq:gamma}
S_{s_\gamma}S_iu_\la\in\Uf S_iu_\la.
\eneq
If $\lan h_i,\lam\ran\ge0$, we have
\eqn
&&\U^-u_\lam\supset
\U^-S_iu_\lam\supset
\U^-S_{s_\gamma}S_iu_\lam
=
S_iS_{s_\beta}u_\lam
\eneqn
Since
$\U^-u_\lam$ is an $\U_i$-module,
it contains $S_{s_\beta}u_\lam$.

Now assume that 
$\lan h_i,\lam\ran<0$.
Then $\lan h_i,s_\beta\lam\ran
=\lan h_i,\lam\ran-\lan\beta^\vee,\la\ran\lan h_i,\beta\ran<0$.
By \eqref{eq:gamma}, we have
\eqn
S_iS_{s_\beta}S_i(u_\infty\otimes t_{s_i\lam}\otimes u_{-\infty})
&=&S_{s_\gamma}(u_\infty\otimes t_{s_i\lam}\otimes u_{-\infty})\\
&&\phantom{88888}\in B(\infty)\otimes t_{s_i\lam}\otimes u_{-\infty}.
\eneqn
Hence applying $S_i^*$ we have
\eqn
&&S_iS_{s_\beta}(u_\infty\otimes t_{\lam}\otimes u_{-\infty})
\in S_i^*(B(\infty)\otimes t_{s_i\lam}\otimes u_{-\infty}),
\eneqn
or equivalently (here $\te^*_i{}^\max b=\tei^*{}^{\eps_i^*(b)}b$ and
$\tfi{}^\max b=\tfi^{\phi_i(b)}b$)
\eqn
S_{s_\beta}(u_\infty\otimes t_{\lam}\otimes u_{-\infty})
&\in& S_iS_i^*\Big(B(\infty)\otimes t_{s_i\lam}\otimes u_{-\infty}\Big)\\
&=&\tf_i{}^\max\te^*_i{}^\max\Big(
B(\infty)\otimes t_{s_i\lam}\otimes u_{-\infty}\Big)\\
&=&\tf_i{}^\max\Big(\bigcup_{n\ge0}
B(\infty)\otimes t_{\lam}\otimes \tei^nu_{-\infty}\Big)\\
&\subset&
B(\infty)\otimes t_{\lam}\otimes u_{-\infty}.
\eneqn
The last inclusion follows from
\eqn
&&\tf_i^\max(b_1\otimes t_\la\otimes b_2)
=b_1'\otimes t_\la\otimes \tf_i^\max b_2
\quad\mbox{for some $b_1'\in B(\infty)$.}
\eneqn
\qed

\subsection{Affine case}
Until now, we have assumed that $\g$ is a symmetrizable Kac-Moody algebra.
From now on, we assume further 
that $\U$ is a quantized {\em affine} algebra.

\subsubsection{Extended Weyl groups}
We take a weight lattice $P$ of rank $\rk(\g)+1$
and an inner product on $P$ as in \cite{KF}.
We set $\td=\Q\otimes P$, which is canonically determined 
by the Dynkin diagram.

Let us define
$\delta\in\sum_i\Z_{\ge0}\al_i$ and $c\in\sum_i\Z_{\ge0}h_i$
by
\eq
\ba{l}
\set{\la\in\sum_i\Z\al_i}{\mbox{$\lan h_i,\la\ran=0$ for every $i\in I$}}
=\Z\delta,\\[5pt]
\set{h\in\sum_i\Z h_i}{\mbox{$\lan h,\al_i\ran=0$ for every $i\in I$}}
=\Z c.
\ea
\eneq
By the inner product of $\td$, we sometimes identify
$\td$ and its dual.
Note that the inner product on $\td$ is so normalized that
$\delta$ and $c$ correspond by this identification.

For $\al\in\Delta^\re$, we set
$c_\al\seteq\max(1,(\al,\al)/2)\in\Z$.
Then we have
$$(\al+\Z\delta)\cap\Delta=\al+c_\al\Z.$$

Let us denote by $P_\cls$ the quotient space $P/(P\cap\Q\,\delta)$,
and let us denote by $\cls\cl P\to P_\cls$ the canonical projection.
Let us denote by $P_\cls^*$ the dual lattice of $P_\cls$, i.e.\ 
$P_\cls^*=\Ker(\delta\cl P^*\to\Z)=(\sum_i\Q\, h_i)\cap P^*$.

Similarly to $P_\cls$, we define
$\td_\cls\seteq\td/\Q\,\delta$, and
let $\cls\cl\td\to\td_\cls$ be the canonical projection.
Define
$\td{}^0\seteq \Ker(c\cl \td\to\Q)$,
and $\td_\cls{}^0=\cls(\td{}^0)$.
The dimension of $\td_\cls{}^0$ is equal to $\rk(\Gg)-1$.
The inner product of $\td$ induces 
a positive definite inner product on $\td_\cls{}^0$.

Let us denote by $\O(\td)$ the orthogonal group,
and $\O(\td)_\delta\seteq\set{g\in\O(\td)}{g\delta=\delta}$ 
the isotropy subgroup at $\delta$.
Then there is an exact sequence

$$1\To \td_\cls{}^0\To[\ t\ ] \O(\td)_\delta\To[\ \cls_0\ ] \O(\td_\cls{}^0)\To 1.$$
Here $t\cl \Gt_\cls{}^0\to \O(\td)_\delta$ is given by

$$\mbox{$t(\cls(\xi))(\la)=\la+(\la,\delta)\xi-(\la,\xi)\delta
-\dfrac{(\xi,\xi)}{2}(\la,\delta)\delta$\quad for $\xi\in\td{}^0$ and $\la\in\td$.}$$
Let us set $W_\cls=\cls_0(W)$.
Then $W_\cls$ is the Weyl group of the root system
$\Delta_\cls\seteq\cls(\Delta^\re)\subset\td_\cls{}^0$.
We define the extended Weyl group $\tW$ by
$$\tW\seteq\set{w\in\O(\td)_\delta}{\mbox{$w\Delta=\Delta$ and
$\cls_0(w)\in W_\cls$}}.$$
Then we have a commutative diagram with the exact rows:
$$
\xymatrix{
1\ar[r] &{\tQ}\ar[r]\ar@{_{(}->}[d] &{W}\ar[r]\ar@{_{(}->}[d]
 &{W_\cls}\ar[r]\ar@{=}[d] &1\\
1\ar[r] &{\tP}\ar@{_{(}->}[d]\ar[r] &{\tW}\ar[r]\ar@{_{(}->}[d] &{W_\cls}\ar[r]
\ar@{_{(}->}[d] &1\\
1\ar[r] &{\td_\cls{}^0}\ar[r]^(.4)t &{\O(\td)_\delta}\ar[r]^{\cls_0}
&{\O(\td_\cls{}^0)}\ar[r] &1
}
$$
Here $\tP$ and $\tQ$ are given by
$$\tP=P_\cls^0\cap P^\vee_\cls{}^0\quad\mbox{and}\quad
\tQ=Q_\cls\cap Q_\cls^\vee,$$
where
\eqn
P_\cls^0&\seteq&\set{\la\in\td_\cls{}^0}
{\mbox{$\lan h_i,\la\ran\in\Z$ for every $i\in I$}},\\
P_\cls^0{}^\vee&\seteq&\set{\la\in\td_\cls{}^0}{
\mbox{$(\al_i,\la)\in\Z$ for every $i\in I$}},\\
Q_\cls&\seteq&\sum_{i\in I}\Z\,\cls(\al_i),\\
Q_\cls^\vee&\seteq&\sum_{i\in I}\Z\,\cls(h_i).
\eneqn
The Weyl group $W$ is a normal subgroup of $\tW$, and 
$\tW$ is a semi-direct product of $W$ and 
$\Aut_0(\Dyn)\seteq\set{\iota}%
{\mbox{$\iota$ is a Dynkin diagram automorphism such that
$\cls_0(\iota)\in W_\cls$}}$.

$$\tP/\tQ\isoto \tW/W\isoto \Aut_0(\Dyn).$$

\begin{remark}
\bi
\item
If $\Gg$ is untwisted, then $(\al,\al)/2\le1$ for every $\al\in\Delta^\re$ and
\eqn
&&\tP=P_\cls^0{}^\vee\subset P_\cls^0{},\quad
\tQ=Q_\cls^\vee\subset Q_\cls.
\eneqn
\item
If $\Gg$ is the dual of an untwisted affine algebra, then
$(\al,\al)/2\ge1$ for every $\al\in\Delta^\re$ and
\eqn
&&\tP=P_\cls^0{}\subset P_\cls^0{}^\vee,\quad
\tQ=Q_\cls\subset Q_\cls^\vee.
\eneqn
\item
If  $\Gg=A^{(2)}_{2n}$, then we have
$(\alpha,\alpha)/2=1/2$, $1$ or $2$, and 
\eqn
&&\tP=\tQ=P_\cls^0=P_\cls^0{}^\vee=Q_\cls=Q_\cls^\vee
=\sum_{\al\in\Delta^\re}\Z\,\cls(\al)
=\sum_{\al\in\Delta^\re,\,(\al,\al)/2=1}\Z\,\cls(\al).
\eneqn
\ei

\end{remark}

\subsubsection{Peter-Weyl theorem}

Let us recall some of the results
by Nakajima and Beck-Nakajima.

The following theorem is conjectured in \cite{KF} by the author
and proved in \cite{BN} by Beck-Nakajima.

\Th[a version of Peter-Weyl theorem]
$$B(\tU)\simeq\Bigl(\bigsqcup_{\la\in P}B(\la)\times B_0(-\la)\Bigr)/W.$$
\enth
Here $B_0(\la)$ is the connected component of $B(\la)$
containing $u_\la$.
Note that $B_0(\la)=B(\la)$
when the level of $\la$ does not vanish.
The Weyl group acts on $\bigsqcup_{\la\in P}B(\la)\times B_0(-\la)$
by $W\ni w\cl B(\la)\times B_0(-\la)
\to B(w\la)\times B_0(-w\la)$
via the action given in \eqref{eq:weyl}.
The left crystal structure
$(\tei,\tfi)$ on $B(\tU)$ is
compatible with the crystal structure of $B(\la)$,
the first factor of $B(\la)\times B_0(-\la)$,
and the right crystal structure
$(\te^*_i,\tf^*_i)$ on $B(\tU)$ is
compatible with the crystal structure of $B_0(-\la)$,
the second factor of $B(\la)\times B_0(-\la)$.

\medskip
For $\la\in P$, there exists a unique symmetric bilinear form
$(\cdot,\cdot)$ on $V(\la)$ that satisfies:
\eqn
&&\mbox{$(u_\la, G(b))=\delta_{b,u_\la}$ for every $b\in B(\la)$,}\\
&&\mbox{$(e_i u,v)=(u,f_i v)$ for every $u$, $v\in V(\la)$,}\\
&&\mbox{$(q(h)u,v)=(u,q(h) v)$ 
for every $u$, $v\in V(\la)$ and $h\in d^{-1} P^*$.}\\
\eneqn

The following theorem is trivial for non-zero level case,
and proved in \cite{N,BN} 
by Nakajima and Beck-Nakajima for the zero level case.

\Th\label{th:nondeg}
\bi
\item
This symmetric bilinear form on $V(\la)$ is non-degenerate.
\item
$(G(b),G(b'))\in q^{(\mu,\mu)-(\la,\la)}(\delta_{b,b'}+\qs\A_0)$ 
for any $\mu\in P$ and
$b,b'\in B(\la)_\mu$.
\item
For $b,b'\in B(\la)_{\la}$, we have
$(G(b),G(b'))=\delta_{b,b'}$.
\ei
\enth
In particular if $v$ is a non-zero vector of $V(\la)$,
then there exists $P\in \U$ such that
$(u_\la,Pv)$ does not vanish.
Note that $(u_\la,Pv)$ coincides with the coefficient of $u_\la$ when we write
$Pv$ as a linear combination of the global basis.

\begin{conjecture}
Theorem \ref{th:nondeg} holds for an arbitrary symmetrizable Kac-Moody 
Lie algebra $\g$.
\end{conjecture}

\section{Extremal vectors}\label{sect:ext}
We assume that $\U$ is a quantized affine algebra.
Let $M$ be an integrable $\U$-module with a global basis.
Let $N$ be a $\Ue$-submodule of $M$ compatible with the global basis of $M$.
Then, for $\la\in P^+$, 
$u_\la\otimes N$ is also a $\Ue$-submodule of $V(\la)\otimes M$
compatible with the global basis.
Hence
$\U(u_\la\otimes N)$ is a $\U$-module compatible with the global basis of
$V(\la)\otimes M$.
\Prop
Assume that for any $b\in B(N)$
if $\eps_i(b)>0$, then $\tfi b\in B(N)\cup\{0\}$.
Or equivalently, for any $i$-string $S$ of $B(M)$,
$S\cap B(N)$ is either $S$ itself, the empty set or
the set consisting of the highest weight vector of $S$.
Assume further that
$f_i N\subset N$ whenever $\lan h_i,\la\ran=0$.
Then we have
$$\U(u_\la\otimes N)\cap u_\la\otimes M=u_\la\otimes N.$$
\enprop
\proof
It is enough to show that
$$B(\U(u_\la\otimes N))\cap u_\la\otimes B(M)=u_\la\otimes B(N).$$
We have
$$B(\U(u_\la\otimes N))
=\set{\tf_{i_1}\cdots\tf_{i_m}b}%
{m\ge 0, i_1,\ldots i_m\in I, b\in u_\la\otimes B(N)}\setminus\{0\}.$$
Hence it is enough to show that, for $b\in B(N)$ such that
$\tfi(u_\la\otimes b)=u_\la\otimes \tfi b$, we have $\tfi b\in B(N)\cup\{0\}$.
Since $\tfi(u_\la\otimes b)=u_\la\otimes \tfi b$
if and only if $\lan h_i,\la\ran\le\eps_i(b)$,
and $\tfi B(N)\subset B(N)\cup\{0\}$ if $\lan h_i,\la\ran=0$,
the assertion follows.
\qed

\Cor\label{cor:ue}
For $\la\in P^+$ and $\mu\in P$ such that
$\lan h_i,\mu\ran\le 0$ whenever $\lan h_i,\la\ran=0$,
we have
$$\U(u_\la\otimes u_\mu)\cap u_\la\otimes V(\mu)=u_\la\otimes\Ue u_\mu,$$
and
$$B\bigl(\U(u_\la\otimes u_\mu)\bigr)
\cap u_\la\otimes B(\mu)=u_\la\otimes B^+(\mu).$$
\encor

Indeed, Lemma \ref{lem:str} allows us to apply the proposition above to
$M=V(\mu)$ and $N=\Ue u_\mu$.

\bigskip
Let $\la\in P^+$ be a dominant integral weight 
and $\mu\in P$ an integral weight.
Then we have a chain of morphisms compatible with global bases
$$\U a_{\la+\mu}\to V(\la)\otimes \U a_\mu\to V(\la)\otimes V(\mu).$$

\Th\label{th:main}
Let $\la\in P^+$ be a dominant integral weight 
and $\mu\in P$ an integral weight.
Then there exist a unique homomorphism
$\U (u_\la\otimes u_\mu)\to V(\la+\mu)$
that sends $u_\la\otimes u_\mu$ to $u_{\la+\mu}$.
Moreover this morphism is compatible with global bases.
\enth
\proof
We have a commutative diagram
$$
\xymatrix{
{\U a_{\la+\mu}}\ar@{>>}[r]\ar@{>>}[dr]
&{\U (u_\la\otimes u_\mu)}\ar@{^{(}->}[r]\ar@{..>}[d]
&{V(\la)\otimes V(\mu)}\\
&{V(\la+\mu)}
}
$$
All the solid arrows are compatible with global bases.
Hence, in order to show the theorem,
it is enough to show the existence of the dotted arrow.

Correspondingly, we obtain the following diagram of crystal bases.
$$
\xymatrix{
{B(\U a_{\la+\mu})}
&{B(\bigl(\U (u_\la\otimes u_\mu)\bigr)}
\ar@{_{(}->}[l]\\
&{B(\la+\mu)}\ar@{_{(}->}[ul]
}
$$
Let $G(b)\in \U a_{\la+\mu}$ be the global basis vector
corresponding to $b\in B(\U a_{\la+\mu})$.
Then $b\in B(\la+\mu)$ means that
$G(b)u_{\la+\mu}\not=0$
(we regard $V(\la+\mu)$ as an $\tU$-module),
and $b\in B\bigl(\U (u_\la\otimes u_\mu)\bigr)$
means $G(b)(u_\la\otimes u_\mu)\not=0$.

Hence we have reduced the problem to the following proposition:
$$\mbox{$B(\la+\mu)\subset B\bigl(\U (u_\la\otimes u_\mu)\bigr)$
as subsets of $B(\U a_{\la+\mu})$.}
$$
Taking $b\in B(\la+\mu)\subset B(\U a_{\la+\mu})$,
let us show that
$b\in B\bigl(\U (u_\la\otimes u_\mu)\bigr)$.

Since $G(b)u_{\la+\mu}\not=0$, Theorem \ref{th:nondeg} implies that
there exists $P\in \U$ such that,
when we write $PG(b)u_{\la+\mu}$ as a linear combinations of
the global basis of $V(\la+\mu)$,
the coefficient of $u_{\la+\mu}$ does not vanish.
Hence, if we write $PG(b)$ as a linear combinations of
the global basis of $\U a_{\la+\mu}$ 
the coefficient of $a_{\la+\mu}$ does not vanish.
Hence $PG(b)(u_\la\otimes u_\mu)$ as a linear combinations of
the global basis of $V(\la)\otimes V(\mu)$
the coefficient of $u_\la\otimes u_\mu$ does not vanish.
We conclude then that $PG(b)(u_\la\otimes u_\mu)\not=0$.
Hence $b\in B\bigl(\U (u_\la\otimes u_\mu)\bigr)$.
\qed

\begin{remark}
Theorem \ref{th:main} holds also for any finite-dimensional $\g$,
because Theorem \ref{th:nondeg} trivially holds in such a case.
\end{remark}

\Cor\label{cor:lmex}
If $\la\in P^+$ and $\mu\in P$, 
then $u_\la\otimes u_\mu\in B(\la)\otimes B(\mu)$
is an extremal vector, and we have an inclusion
$B(\la+\mu)\subset B(\la)\otimes B(\mu)$
as subsets of $B(\la)\otimes B(\U a_\mu)$.
\encor

Note that the first statement 
holds for an arbitrary Kac-Moody Lie algebra.

\section{Fundamental representations}
Write the smallest positive imaginary root $\delta$
and the smallest positive imaginary coroot $c$ as
\eqn
&&\delta=\sum_ia_i\alpha_i,\qquad
c=\sum_ia_i^\vee h_i.
\eneqn
Then we have
$$a^\vee_i=\dfrac{(\al_i,\al_i)}{2}a_i.$$

We choose $0\in I$ such that 
\bi
\item
Setting $I_0=\set{i\in I}{i\not=0}$ and
$W_0\seteq \lan s_i\mathbin;i\in I_0\ran\subset W$,
the composition $W_0\hookrightarrow W\To[\cls_0]W_\cls$ is an
isomorphism.
\item
$a_0=1$.
\ei

Such a $0$ exists and is unique up to a Dynkin diagram automorphism.

If $\g=A^{(2)}_{2n}$,
$\alpha_0$ is the longest simple root and $a^\vee_0=2$.

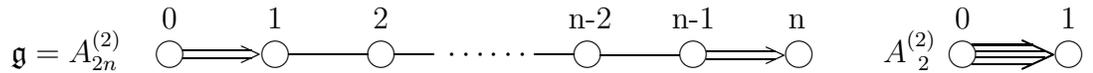
\begin{figure}[H]
\begin{center}
\begin{picture}(100,110)(90,-80)

\put(-60,-3){\mbox{$\g=A^{(2)}_{2n}$}}
\put(-3,10){0}
\put(37,10){1}
\put(77,10){2}

\put(0,0){\circle{10}}
\put(40,0){\circle{10}}
\put(80,0){\circle{10}}
\put(5,-1.5){\line(1,0){26}}
\put(5,1.5){\line(1,0){26}}
\put(26,-3){$>$}
\put(44.8,0){\line(1,0){30}}
\put(85,0){\line(1,0){15}}
\put(105,-3){$\cdots\cdots$}
\put(3,0){
\put(148,10){n-2}
\put(187,10){n-1}
\put(231,10){n}

\put(155,0){\circle{10}}
\put(195,0){\circle{10}}
\put(235,0){\circle{10}}

\put(135,0){\line(1,0){15}}
\put(160,-.3){\line(1,0){30}}
\put(200,-1.8){\line(1,0){26}}
\put(200,1.2){\line(1,0){26}}
\put(221.5,-3.3){$>$}}

\put(300,0){
\put(-30,-3){\mbox{$A^{(2)}_{\,\,\,2}$}}
\put(0,0){\circle{10}}
\put(40,0){\circle{10}}
\put(3.5,-3.6){\line(1,0){23}}
\put(5,-1.15){\line(1,0){26}}
\put(5,1.15){\line(1,0){26}}
\put(3.5,3.6){\line(1,0){23}}
\put(20.5,-5.2){{\LARGE\mbox{$>$}}}
\put(-3,10){0}
\put(37,10){1}
}
\put(80,20){
\put(0,-50){\mbox{$\delta=\al_0+2(\al_1+\cdots+\al_n)$,}}
\put(0,-70){\mbox{$c=2(h_0+\cdots+h_{n-1})+h_n$.}}
\put(-50,-90){\mbox{$(\al_0,\al_0)=4$, 
$(\al_n,\al_n)=1$, $(\al_i,\al_i)=2$ for $0<i<n$}}
}
\end{picture}
\end{center}
\caption{$\g=A^{(2)}_{2n}$}
\end{figure}
If $\g$ is not of type $A^{(2)}_{2n}$, then $a^\vee_0=1$.
Note that $\delta-\alpha_0\in \Delta^\re$ if $\g$ is not of type $A^{(2)}_{2n}$,
and $(\delta-\alpha_0)/2\in \Delta^\re$ if $\g$ is of type $A^{(2)}_{2n}$.
Hence one has always $s_{\delta-\alpha_0}\in W$.

Let us denote by $U_q(\g_0)$ the subalgebra of $\U$ generated by $e_i$, $f_i$
($i\in I_0$). This is the quantized universal enveloping algebra
associated with a finite-dimensional simple Lie algebra.

Let $\vp_k$ be a fundamental weight of level $0$. 
That is, $\{\vp_k\}_{k\in I_0}$
is a set of vectors such that $\lan h_j,\vp_k\ran=0$ for 
$j\in I_0$ with $j\not=k$, and

\eqn
P_\cls^0{}^+&\seteq&
\set{\la\in \Gt^*_{\cls}}{\mbox{$\lan c,\la \ran=0$ and 
$\lan h_i,\la\ran\in\Z_{\ge0}$ for every $i\in I_0$}}\\
&=&\sum_{k\in I_0}\Z_{\ge0}\,\vp_k.
\eneqn

A fundamental weight of level $0$ is unique up to 
$\Q\,\delta$.
We can take
\eq
\vp_k&=&\begin{cases}
\La_k-a_k^\vee\La_0&\mbox{when $a_0^\vee=1$,}\\
\dfrac{2}{(\al_k,\al_k)}\La_k-\La_0&\mbox{when $\g=A^{(2)}_{2n}$.}
\end{cases}\eneq
Here $\Lambda_k$ is a vector in $P$ satisfying
$\lan h_i,\Lambda_k\ran=\delta_{ik}$ for $i\in I$.

\smallskip
Let $k\in I\setminus\{0\}$.
Set $c_k=\max(1,(\alpha_k,\alpha_k)/2)\in\Z$.
Then we have
$$\set{n\in\Z}{\alpha_k+n\delta\in\Delta}=\Z c_k,$$
and 
$$W\vp_k\cap (\vp_k+\Z\delta)=\vp_k+\Z c_k\delta.$$
We have
$\vp_k+n\delta\in \Wt(V(\vp_k))$ if and only if $n\in c_k\Z$.

For any $\xi=w\la\in W\la$, we use the notation
$u_\xi$ for the extremal vector $\Sn_w u_\la\in V(\la)$.
Note that $\Sn_w u_\la$ is a unique global basis vector of weight $\xi$.

We denote by $\Us$ the subalgebra of $\U$ generated by $e_i$, $f_i$ ($i\in I$)
and $q(h)$ ($h\in d^{-1}P_\cls^*\subset d^{-1}P^*$).


Then there exists a unique $\Us$-morphism
$z_k\cl V(\vp_k)\to V(\vp_k)$ sending
$u_{\vp_k}$ to $u_{\vp_k+c_k\delta}$.
The operator $z_k$ has weight $c_k\delta$.
The global basis of $V(\vp_k)$ is stable by $z_k$.
We have
$$\mbox{$z_ku_\xi=u_{\xi+c_k\delta}$\quad for every $\xi\in W\cdot\la$.}$$

The quotient $W(\vp_k)\seteq V(\vp_k)/(z_k-1)V(\vp_k)$
is an irreducible $\Us$-module with a global basis.
The morphism $V(\vp_k)\epi W(\vp_k)$ sends the members of the global basis of
$V(\vp_k)$ to the one of $W(\vp_k)$.

In this section, we set
\eq
\la\seteq\vp_k,\quad 
\mu\seteq w_0\vp_k\quad\mbox{where $w_0$ be the longest element of $W_0$.}
\eneq
Then $\mu\equiv -\vp_{k'}\mod \Z\delta$ for some $k'\in I_0$.

Then $U_q(\g_0)u_\la=U_q(\g_0)u_\mu$ is an irreducible
$U_q(\g_0)$-module with highest weight $\la$ and lowest weight $\mu$.
Note that $\Ue u_\mu$ is a $U_q(\g_0)$-module.
We have
$$\mbox{$u_\xi\in\Ue u_\mu$ for any 
$\xi\in W_0\cdot\la=W\cdot\la\cap (\la+\sum\limits_{i\in I_0}\Z\alpha_i)$.}
$$

\Lemma
$z_ku_\mu\in\Ue u_\mu$.
\enlemma
\proof
Assume $a_0^\vee=1$.
Then we have
$\tilde \alpha_{k'}\seteq c_k\alpha_{k'}^\vee\in \tQ$
and $t(\tilde\al_{k'})(\mu)=\mu+c_k\delta$.
Hence $z_ku_\mu=S_{t(\tilde\al_{k'})}^\norm u_\mu$.
Since we have $c_k\delta-\al_{k'}\in\Delta^\re_+$,
$t(\tilde\al_{k'})=s_{c_k\delta-\al_{k'}}s_{\al_{k'}}$,
and $(\al_{k'},\mu)=(c_k\delta-\al_{k'},s_{k'}\mu)<0$,
Proposition \ref{prop:beta} 
implies that $S_{t(\tilde\al_{k'})}^\norm u_\mu\in\Ue u_\mu$.

Now assume that $\g=A^{(2)}_{2n}$.
Then $t(-\al_0/2)\mu=\mu+\delta$ and 
$t(-\al_{0}/2)=s_{\al_{0}}s_{(\delta-\al_{0})/2}$,
and $(\delta-\al_{0},\mu)=(\al_{0},s_{(\delta-\al_0)/2}\mu)<0$.
Note that $(\delta-\al_0)/2\in\Delta^\re_+$.
Hence Proposition \ref{prop:beta}
implies that $S_{t(\tilde\al_{0})}^\norm u_\mu\in\Ue u_\mu$.
\qed

\Lemma\label{lem:ueW}
\bi
\item
$B^+(\mu)\setminus z_kB^+(\mu)\simeq B(W(\la))$
as a crystal over $\g_0$.
\item
$\Ue u_\mu/\Ue z_ku_\mu$ is isomorphic to $W(\la)$ as a $U_q(\g_0)$-module
\ei
\enlemma
\proof
The crystal $B^+(\mu)$ is a regular crystal
over $\g_0$.
The crystal $B^+(\mu)$ is invariant by $z_k$, and 
$\bigcup\limits_{n\in\Z}z_k^n B^+(\mu)=B(\mu)$, 
$\bigcap\limits_{n\in\Z}z_k^n B^+(\mu)=\emptyset$.
On the other hand, the crystal $B(W(\la))$ is isomorphic to the quotient of
$B(\mu)$ by the action of $\Z$ given by $z_k$.
Hence $B^+(\mu)\setminus z_kB^+(\mu)\to B(W(\la))$
is bijective.
(ii)  follows from (i).
\qed

\Lemma\label{lem:4.2}
Assume that $\lan h_0,\la\ran=-1$ 
\ro i.e.\ $a^\vee_k=1$ or $\g=A^{(2)}_{2n}$\rf.
\bi
\item
$c_k=1$.
\item $W(\la)$ is an irreducible $U_q(\g_0)$-module.
\ei
\enlemma
\proof
(i)\quad if $\g=A^{(2)}_{2n}$, then $(\al_k,\al_k)/2\le1$ for $k\not=0$,
and hence $c_k=1$. If $a^\vee_k=1$, then
$1=a^\vee_k=\dfrac{(\alpha_k,\al_k)}{2}a_k\ge \dfrac{(\alpha_k,\al_k)}{2}$,
and hence $c_k=1$.

\noindent 
(ii)\quad By Lemma \ref{lem:ueW}, 
it is enough to show that $N:=\Ue u_\mu/\Ue z_k u_\mu$ is irreducble as
a $U_q(\g_0)$-module.
Note that $e_0u_\la=u_{s_o\la}=u_{\la+\delta}=z_ku_\la\in \Ue z_ku_\mu$
and $e_i u_\la=0$ for $i\not=0$.
Since $\Ue u_\mu=\Ue U_q(\g_0) u_\la=U_q(\g_0)\Ue u_\la$,
we have
$N=U_q(\g_0)u_\la \mod z_k\Ue u_\mu$.
\qed

\medskip
Two vectors $u_\L\otimes u_\la$
and $u_\L\otimes u_\mu$ are extremal vectors
in the same connected component of $B(\L)\otimes B(\la)$.
Since the level of $\L+\mu$ is equal to one, 
there exists a unique dominant weight $\xi_0$ of level one such that
$\xi_0\in W\cdot(\L+\la)=W\cdot(\L+\mu)$.
By Corollary \ref{cor:lmex}, the connected component
of $B(\L)\otimes B(\la)$ containing $u_\L\otimes u_\la$
is isomorphic to $B(\xi_0)$.
Set $M=V(\L)\otimes V(\mu)$
and $M_n=\U(u_\L\otimes z_k^nu_{\la})=
\U(u_\L\otimes z_k^nu_{\mu})\subset M$.
Note that Corollary \ref{cor:ue} implies 
$$M_n\cap \bigl(u_\L\otimes V(\la)\bigr)=u_\L\otimes \Ue z_k^nu_{\mu}.$$

\Lemma 
\bi
\item
$M_{n+1}\subset M_n$.
\item
$M=\bigcup_{n\in\Z} M_n$.
\item $\bigcap_{n\in \Z}M_n=0$.
\ei
\enlemma
\proof
(i) is obvious.
Since $W(\mu)$ is generated by $u_\mu$ as an $\Ue$-module 
(by \cite[Proposition 1.16]{AK}),
we have $V(\mu)=\cup_n\Ue z_k^nu_{\mu}$, which implies (ii).

In order to prove (iii), it is enough to show that
$\cap_n B(M_n)=\emptyset$.
Any vector $b\in \cap_n B(M_n)$ is connected with a vector in 
$u_\L\otimes B(\mu)$.
Since $ B(M_n)\cap \bl(u_\L\otimes B(\mu)\br)=u_\L\otimes z_k^nB^+(\mu)$,
the result follows from
$\cap_nz_k^nB^+(\mu)=\emptyset$,
which is an immediate consequence of
$\cap_n(\mu+n\delta+Q_+)=\emptyset$.
\qed

\Lemma\label{lem:MW}
$B(M_0)\setminus B(M_1)\simeq B(\L)\otimes B(W(\la))$.
\enlemma
 The proof is similar to the one of Lemma \ref{lem:ueW}.

\Prop
\bi
\item
The vector $u_\L\otimes u_{\la+nc_k\delta}$,
regarded as a vector of $M_n/M_{n+1}$,
is an extremal vector.
\item
\eqn
\xi_0&\equiv&
\left\{
\ba{lcl}
\Lambda_k&&\mbox{if $\lan h_0,\la\ran=-1$,
i.e.\ $a^\vee_k=1$ or $\g=A^{(2)}_{2n}$,}\\[10pt]
\iota^{-1}(\L)&&\parbox{310pt}{if $a^\vee_0=1$,
$\dfrac{(\al_k,\al_k)}{2}\ge1$, 
and $\iota$ is the Dynkin diagram automorphism 
such that $t(\la)\in W\iota$,}\\[20pt]
\Lambda_4&&\mbox{if $\g=F^{(1)}_4$ and $k=3$,}
\ea\right.
\eneqn
$\mod \Q\,\delta$.
For the last case, see {\rm Figure \ref{fig:F}} in the proof.
\ei
\enprop

\proof
We may assume that $n=0$.
If $u_\L\otimes u_{\la}\in M_1$,
then $u_{\la}\in\Ue z_ku_\mu$, which is a contradiction.
Hence $u_\L\otimes u_{\la}\mod M_1$ is a non-zero vector.

We divide the proof into three cases.

\smallskip
\noindent
Case 1)\quad $a^\vee_k=1$ or $\g=A^{(2)}_{2n}$

In this case $\lan h_0,\la\ran=-1$
and $c_k=1$ by Lemma \ref{lem:4.2}.
We shall show that $u_\L\otimes u_\la$ is a highest weight vector
of $M_0/M_{1}$
We have
$e_i(u_\L\otimes u_\la)=0$ for $i\not=0$.
We have
$e_0(u_\L\otimes u_\la)=u_\L\otimes e_0u_\la$,
and $e_0u_\la=S_0u_\la=u_{s_{\al_0-\delta}\la+\delta}
=z_ku_{s_{\al_0-\delta}\la}\in z_k \Ue u_\mu$.

\smallskip
\noindent
Case 2)\quad $a^\vee_0=1$ and $\dfrac{(\al_k,\al_k)}{2}\ge1$

We have $\la\in\tP$ and hence $t(\la)\in \tW$.
We have then
$t(-\la)(\L+\la)\equiv\La_0 \mod\Z\delta$.
Hence if we write $t(\la)=w\iota$ with $w\in W$ and 
a Dynkin diagram automorphism $\iota$, then
we have 
$\xi_0=w^{-1}(\L+\la)=\iota^{-1}t(-\la)(\L+\la)\equiv \iota^{-1}(\L)
\mod\Q\,\delta$.

Set $v=S_{t(\la)}^{-1}(u_\L\otimes u_\la)$.
Here we regard $S_{\iota^{-1}}$ as an isomorphism
$V(\L)\otimes V(\la)\to V(\iota^{-1}\L)\otimes V(\iota^{-1}\la)$
such that $S_{\iota^{-1}}(au)=\iota^{-1}(a)S_{\iota^{-1}}(u)$ 
for $a\in \U$ and $u\in V(\L)\otimes V(\la)$.
Hence, $v$ is regarded as a vector in
$V(\iota^{-1}\L)\otimes V(\iota^{-1}\la)$.
We shall show that
$e_i v\in S_{\iota^{-1}}M_1$ for every $i\in I$.

\bi
\item
$i\not=0,k$

In this case $t(\la)\al_i=\al_i$ holds.
Hence we have $T_{t(\la)}e_i=e_i$, and
$$S_{t(\la)}e_iv=(T_{t(\la)}e_i)(u_\L\otimes u_\la)=e_i(u_\L\otimes u_\la)=0.$$

\item $i=k$

\noindent
Since $\lan h_k,\wt(v)\ran=0$,
it is enough to show that
$e_k S_k^{-1}S_{t(\la)}^{-1}(u_\L\otimes u_\la)\in S_{\iota^{-1}}M_1$.
This is equivalent to saying that
$(T_{t(\la)s_k}e_k)(u_\L\otimes u_\la)\in M_1$.
Since $t(\la)s_k\alpha_k=c_k\delta-\alpha_k\in \Delta^+$,
we have $T_{t(\la)s_k}e_k\in \Ue$ and
$(T_{t(\la)s_k}e_k)(u_\L\otimes u_\la)
=u_\L\otimes (T_{t(\la)s_k}e_k)u_\la$.
The last factor is calculated as
$(T_{t(\la)s_k}e_k)u_\la=S_{t(\la)s_k}e_k S_{t(\la)s_k}^{-1}u_\la
=S_{t(\la)s_k}e_k u_{s_kt(-\la)\la}$ up to a non-zero constant multiple.
Since we have
$\lan h_k, s_kt(-\la)\la\ran=-\lan h_k,\la\ran=-1$,
we obtain $e_k u_{s_kt(-\la)\la}=u_{t(-\la)\la}$.
Thus we obtain $(T_{t(\la)s_k}e_k)u_\la=S_{t(\la)s_k}u_{t(-\la)\la}
=u_{t(\la)s_kt(-\la)\la}$ up to a non-zero constant multiple.
Since $t(\la)s_kt(-\la)\la=s_{t(\la)\alpha_k}\la=s_{\al_k-c_k\delta}\la
=s_k\la+c_k\delta$, we have
$$u_{t(\la)s_kt(-\la)\la}=u_{s_k\la+c_k\delta}=z_ku_{s_k\la}\in z_k\Ue u_\mu.$$

\item $i=0$

Let us first show that
$$\qquad\mbox{$v=S_{t(w\la)}^{-1}(u_\L\otimes u_{w\la})$ 
up to a non-zero constant multiple for every $w\in W_0$,}$$
by the induction of the length of $w$.
Assuming that the assertion is true, we shall show that
it is true for $s_iw$ for $i\in I_0$.
According that $s_it(w\la)\gtrless t(w\la)$, we have
$S_{t(s_iw\la)}=S_i^{\pm}S_{t(w\la)}S_i^\mp$.
Hence we have
\eqn
S_{t(s_iw\la)}^{-1}(u_\L\otimes u_{s_iw\la})
&=&S_i^{\pm}S_{t(w\la)}^{-1}S_i^\mp(u_\L\otimes u_{s_iw\la})\\
&=&S_i^{\pm}S_{t(w\la)}^{-1}(u_\L\otimes u_{w\la})
=S_i^{\pm}v=v
\eneqn
up to a non-zero constant multiple.

Now we divide the proof into two cases.

\begin{enumerate}
\item $(\al_k,\al_k)/2=1$

 In this case, $c_k=1$.
There exists $w\in W_0$ such that
$$w^{-1}\alpha_0\equiv -\al_k \mod \Z\delta.$$
Hence we have $\lan h_0, w\la\ran=-1$,
which implies that
$t(w\la)\al_0=\al_0+\delta$
and $(t(w\la)\al_0,w\la)=-1$.
Then we have
\eqn
S_{t(w\la)}e_0v
&=&S_{t(w\la)}e_0S_{t(w\la)}^{-1}(u_\L\otimes u_{w\la})\\
&=&u_\L\otimes S_{t(w\la)}e_0S_{t(w\la)}^{-1}u_{w\la}\\
&=&u_\L\otimes S_{t(w\la)}S_0S_{t(w\la)}^{-1}u_{w\la}\\
&=&u_\L\otimes S_{\alpha_0+\delta}u_{w\la}\\
&=&u_\L\otimes u_{s_{\alpha_0+\delta}w\la}.
\eneqn
On the other hand, we have
\eqn
s_{\alpha_0+\delta}w\la
&=&w\la-(\alpha_0+\delta,w\la) (\al_0+\delta)\\
&=&w\la-(\alpha_0-\delta,w\la)(\al_0-\delta)+2\delta\\
&=&s_{\delta-\al_0}w\la+2\delta.
\eneqn
This implies that
$$u_{s_{\alpha_0+\delta}w\la}=z_k^2 u_{s_{\delta-\al_0}w\la}
\in z_k\Ue u_\mu.$$
Hence $e_0(u_\L\otimes u_\la)\in S_{t(w\la)}^{-1}M_1=S_{\iota^{-1}}M_1$.

\item $(\al_k,\al_k)/2>1$

 In this case, by the classification of affine Dynkin diagrams,
there exists $i\not=0,k$ such that $(\al_i,\al_i)/2=1$.
Let us take $w\in W_0$ such that
$w^{-1}\al_0\equiv\al_i \mod \Z\delta$.
Hence $(w^{-1}\al_0.\la)=0$, which implies $t(w\la)\al_0=\al_0$.
Hence we have $S_{t(w\la)}e_0S_{t(w\la)}^{-1}=e_0$ and
\eqn
S_{t(w\la)}e_0v&=&
S_{t(w\la)}e_0S_{t(w\la)}^{-1}(u_\L\otimes u_{w\la})\\
&=&u_\L\otimes e_0u_{w\la}.
\eneqn
Since $w\la+\alpha_0\in w(\la+\al_i)+\Z\delta$
is not a weight of $V(\la)$, $e_0u_{w\la}$ must vanish.
\end{enumerate}
\ei

\smallskip
\noindent
Case 3) the remaining case 
(i.e.\ $a_0^\vee=1$, $a_k^\vee>1$ and $(\al_k,\al_k)/2<1$)

By the classification of affine Dynkin diagrams, there is only one
remaining case, namely $\g=F^{(1)}_4$ and $k=3$:

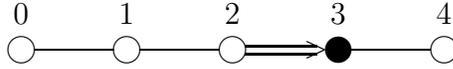
\begin{figure}[H]
\begin{center}
\begin{picture}(100,70)(30,-50)

\put(-3,10){0}
\put(37,10){1}
\put(77,10){2}
\put(117,10){3}
\put(157,10){4}

\put(0,0){\circle{10}}
\put(40,0){\circle{10}}
\put(80,0){\circle{10}}
\put(120,0){\circle*{10}}
\put(160,0){\circle{10}}
\put(4.8,0){\line(1,0){30}}
\put(44.8,0){\line(1,0){30}}
\put(84.8,1.5){\line(1,0){27}}
\put(84.8,-1.5){\line(1,0){27}}
\put(107,-3){$>$}
\put(124.8,0){\line(1,0){30}}

\put(-5,0){
\put(0,-30){
\mbox{$\delta=\al_0+2\al_1+3\al_2+4\al_3+2\al_4$,}
}
\put(2,-47){
\mbox{$c=h_0+2h_1+3h_2+2h_3+h_4$,}
}
}
\end{picture}
\end{center}
\caption{$\g=F^{(1)}_4$}
\label{fig:F}
\end{figure}
We have in this case $c_3=1$, $\la=\La_3-2\La_0$ and
$\mu=\la-4\al_1-8\al_2-12\al_3-6\al_4$.

We have
$$s_4s_3s_2s_1s_0(\L+\la)=\L+\la+\al_0+\al_1+\al_2+\al_4+\al_4\equiv \La_4\mod\Z\delta.$$
Set $x=s_0s_1s_2s_3s_4$ and $v=S_x^{-1}(u_\L\otimes u_\la)$.
Let us show that $v$ is a highest weight vector of $M_0/M_1$,
i.e.\ $e_iv\in M_1$ for $i\in I$.

\begin{enumerate}
\item $i=0$

Since $x\al_0=\al_1$, we have
$S_x e_0v=(T_x e_0)(u_\L\otimes u_\la)=e_1(u_\L\otimes u_\la)=0$.
\item $i=1$

Since $x\al_1=\al_2$, we can conclude $e_1v=0$ by the same argument as above.
\item $i=2$

We have $x\al_2=\al_0+\al_1+\al_2+2\al_3$.
Hence we have
\eqn
S_x e_2v&=&S_x e_2S_x^{-1}(u_\L\otimes u_\la)\\
&=&u_\L\otimes(S_x e_2S_x^{-1}u_\la),
\eneqn
and since $\lan xh_2,\la\ran=-1$, we have
\eqn
S_x e_2S_x^{-1}u_\la&=&S_x S_2S_x^{-1}u_\la
=u_{s_{x\alpha_2}\la}\\
&=&u_{\la+\al_0+\al_1+\al_2+2\al_3}\\
&=&u_{\la+\delta-\al_1-2\al_2-2\al_3-2\al_4}\\
&=&z_ku_{\la-\al_1-2\al_2-2\al_3-2\al_4}\in z_k\Ue u_\mu.
\eneqn
This implies that $e_2v\in M_1$.

\item $i=3$

$e_3v=0$ follows from $x\al_3=\al_4$.
\item $i=4$

We have seen that $v\in M_0/M_1$ is invariant by $S_2$ and $S_3$.
Hence it is enough to show that $e_4 S_3^{-1}S_2^{-1}v=0$.
Since
$xs_2s_3\al_4=\al_3$,
we have
$$S_{xs_2s_3}e_4 S_3^{-1}S_2^{-1}v
=(T_{xs_2s_3}e_4)(u_\L\otimes u_\la)
=e_3 (u_\L\otimes u_\la)=0.$$
\end{enumerate}
\qed

\Th\label{th:main1}
$$\U(u_\L\otimes u_{\la})/
\U(u_\L\otimes u_{\la+c_k\delta})\simeq 
V(\L+\la).$$
By this isomorphism,
$u_\L\otimes u_{\la}$ corresponds to $u_{\L+\la}$.
\enth
\proof
By the preceding proposition, there exists a morphism
$g\cl V(\L+\la)\to M_0/M_1$, 
sending $u_{\L+\la}$ to $u_\L\otimes u_\la\mod M_1$.
On the other hand, Theorem \ref{th:main} implies the existence of
a morphism $\psi\cl M_0\to V(\L+\la)$. 
Since $\L+\la+c_k\delta$ is not a weight of $V(\L+\la)$,
$\psi$ factors through $M_0/M_1$ and thus we obtain a morphism
$M_0/M_1\to V(\L+\la)$ sending $u_\L\otimes u_\la$ to $u_{\L+\la}$.
Obviously it is an inverse of $g$.
\qed

\medskip
Note that the theorem holds if we replace $\la$ with $\mu$.

\Cor\label{cor:fund}
$u_\L\otimes \bigl(B^+(\mu)\setminus B^+(\mu+c_k\delta)\bigr)
\simeq B^+(\L+\mu)$.
In particular 
$ B(W(\la))\simeq B^+(\L+\mu)$
as a crystal over $\g_0$.
\encor
\proof
By the preceding theorem, we have
\eqn
&&\Ue u_{\L+\mu}\simeq
\Ue(u_\L\otimes u_{\mu})/
\bigl(\U(u_\L\otimes u_{\mu+c_k\delta})\cap \Ue(u_\L\otimes u_{\mu})\bigr)
\eneqn
On the other hand,
Corollary \ref{cor:ue} implies that
$\U(u_\L\otimes u_{\mu+c_k\delta})\cap \Ue(u_\L\otimes u_{\mu})
\subset \U(u_\L\otimes u_{\mu+c_k\delta})\cap u_\L\otimes V(\mu)
=u_\L\otimes \Ue u_{\mu+c_k\delta}$, which implies that
$\U(u_\L\otimes u_{\mu+c_k\delta})\cap \Ue(u_\L\otimes u_{\mu})
=u_\L\otimes \Ue u_{\mu+c_k\delta}$.
Hence we have
$$\Ue u_{\L+\mu}\simeq
\bigl(u_\L\otimes \Ue u_{\mu}\bigr)
/\bigl(u_\L\otimes \Ue u_{\mu+c_k\delta}\bigr).$$
Thus we obtain the desired result.
\qed

\medskip
Lemma \ref{lem:MW} and Theorem \ref{th:main1} imply the following result.

\Cor
$B(\L)\otimes B(W(\la))\simeq B(\L+\la)$.
\encor

\Cor
There exists a unique vector $b\in B(W(\la))$ such that
$\eps_i(b)\le\delta_{i,0}$,
\encor
\proof
The condition is equivalent to saying that
$u_\L\otimes b$ is a highest weight vector,
and the preceding corollary implies that
$B(\L)\otimes B(W(\la))$ has a unique highest weight vector.
\qed

\Cor
\bi
\item
If $b\in B(\mu)$ satisfies $\wt(b)\not\in \xi_0-\L-c_k\delta+Q_-$,
then $b\in B^+(\mu)$.
\item If an integral weight $\eta$ satisfies
$\eta\not\in \xi_0-\L-c_k\delta+Q_-$, then
$V(\la)_\eta=(\Ue u_\mu)_\eta$.
\ei
\encor
\proof
There exist $b'\in B^+(\mu)\setminus B(\mu+c_k\delta)$ and $n\in\Z$ such that
$b=z_k^nb'$.
By  Corlollary \ref{cor:fund}, we have
$\wt(u_\L\otimes b')=\L+\wt(b)-nc_k\delta\in\xi_0+Q_-$.
Hence the assumption implies $n\ge0$.
Thus we conclude (i), and (ii) follows from (i).
\qed

\end{document}